\documentclass{article}
\usepackage{amssymb,latexsym,amsmath,amsfonts,amsthm}
\usepackage{graphicx,pstricks}
\newtheorem{theorem}{Theorem}
\newtheorem{lemma}[theorem]{Lemma}

\newcommand{\Ary}{\mathcal A}
\newcommand{\sA}{{\hat A}}
\newcommand{\Bin}{\mathcal B}
\newcommand{\For}{\mathcal F}
\newcommand{\D}{\mathcal D}
\newcommand{\E}{\mathcal E}
\newcommand{\G}{\mathcal G}
\newcommand{\Pla}{\mathcal P}
\newcommand{\PF}{\mathcal{P\!F}}
\newcommand{\Q}{\mathcal Q}
\newcommand{\T}{\mathcal T}
\newcommand{\1}{\bar 1}
\newcommand{\2}{\bar 2}
\newcommand{\one}{\bf b}
\newcommand{\two}{\bf w}
\newcommand{\bi}{\mathrm{bi}}
\newcommand{\bt}{\bar t}
\newcommand{\f}{{\bf f}}
\newcommand{\pv}{\mathrm{pv}}
\newcommand{\ubt}{\mathfrak b}
\newcommand{\lbt}{\ubt,{\omega}}
\newcommand{\lhs}{\mathrm{LHS}}
\newcommand{\rhs}{\mathrm{RHS}}

\begin{document}
\title{A combinatorial proof of Postnikov's identity and a generalized enumeration of labeled trees}
\author{%
Seunghyun Seo\thanks{Department of Mathematics, Brandeis University, Waltham, MA 02454-9110, USA. Email:\,shseo@brandeis.edu. Research supported by the Post-doctoral Fellowship Program of Korea Science and Engineering Foundation (KOSEF).}}
\maketitle
\begin{abstract}
In this paper,  we give a simple combinatorial explanation of a formula of A.~Postnikov relating bicolored rooted trees to bicolored binary trees. We also present generalized formulas for the number of  labeled $k$-ary trees, rooted labeled trees, and labeled plane trees.
\end{abstract}

\section{Introduction}
In Stanley's 60th Birthday Conference, Postnikov~\cite[p.~21]{Po} showed the following identity:
\begin{equation}
(n+1)^{n-1} = \sum_{\ubt} \,\dfrac{n!}{2^n} \prod_{v \in V\!(\ubt)} \left( 1+ \dfrac{1}{h(v)} \right), 
\label{eqn:ori}
\end{equation} 
where the sum is over unlabeled binary trees $\ubt$ on $n$ vertices and $h(v)$ denotes the number of descendants of $v$ (including $v$). The figure below illustrates all~$5$ unlabeled binary trees on $3$ vertices, with the value of $h(v)$ assigned to each vertex $v$. In this case, identity~(\ref{eqn:ori}) says that $(3+1)^2=3+3+4+3+3$\,.
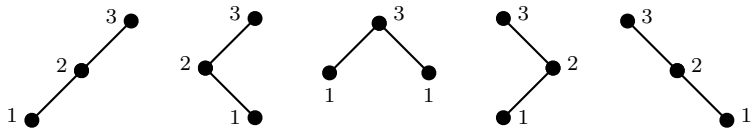
\begin{center}
\begin{pspicture}(0,0)(10,2)\footnotesize
\psset{linewidth=1\pslinewidth}
\psset{dotsize=.2 0} 
\psset{xunit=.66cm, yunit=.66cm}
\rput(2.5,2.45){%
 \psline{*-*}(0,0)(-1,-1)
 \psline{*-*}(-1,-1)(-2,-2)
 \rput(-0.4,0.1){$3$}
 \rput(-1.4,-0.9){$2$}
 \rput(-2.4,-1.9){$1$}
}
\rput(5,2.5){%
 \psline{*-*}(0,0)(-1,-1)
 \psline{*-*}(-1,-1)(0,-2)
 \rput(-0.4,0.1){$3$}
 \rput(-1.4,-0.9){$2$}
 \rput(-0.4,-2){$1$}
}
\rput(7.5,2.4){%
 \psline{*-*}(0,0)(-1,-1)
 \psline{*-*}(0,0)(1,-1)
 \rput(0.4,0.2){$3$}
 \rput(-1,-1.45){$1$}
 \rput(1,-1.45){$1$}
}
\rput(10,2.5){%
 \psline{*-*}(0,0)(1,-1)
 \psline{*-*}(1,-1)(0,-2)
 \rput(0.4,0.1){$3$}
 \rput(1.4,-0.9){$2$}
 \rput(0.4,-2){$1$}
}
\rput(12.5,2.45){%
 \psline{*-*}(0,0)(1,-1)
 \psline{*-*}(1,-1)(2,-2)
 \rput(0.4,0.1){$3$}
 \rput(1.4,-0.9){$2$}
 \rput(2.4,-1.9){$1$}
}
\end{pspicture}
\end{center}
Postnikov derived this identity from the study of  a combinatorial interpretation for mixed Eulerian numbers, which are coefficients of certain reparametrized {\em volume polynomials} which introduced in~\cite{PS}. For more information, see~\cite{PS,Po}.  

In the same talk, he also asked for a combinatorial proof of identity~(\ref{eqn:ori}). Multiplying both sides of~(\ref{eqn:ori}) by $2^n$ and expanding the product in the right-hand side yields 
\begin{equation}
2^{n}\,(n+1)^{n-1} =  \sum_{\ubt} n! \!\!\sum_{\alpha \subseteq V\!(\ubt)} \prod_{v \in \alpha} \dfrac{1}{h(v)}\,.  
\label{eqn:main}
\end{equation}
Let $\lhs_n$ (resp.~$\rhs_n$) denote the left-hand (resp.~right-hand) side of  (\ref{eqn:main}).

The aim of this paper is to find a combinatorial proof of (\ref{eqn:main}). In section~\ref{sec:lhsrhs} we construct the sets $\For_n^{\bi}$ of labeled bicolored forests on~$[n]$ and $\D_n$ of certain labeled bicolored binary trees, where the cardinalities equal $\lhs_n$ and $\rhs_n$, respectively. In section~\ref{sec:bij} we give a bijection between $\For_n^{\bi}$ and $\D_n$, which completes the bijective proof of~(\ref{eqn:main}). Finally, in section~\ref{sec:gen}, we present generalized formulas for the number of labeled $k$-ary trees, rooted labeled trees, and labeled plane trees.      

\section{Combinatorial objects for $\lhs_n$ and $\rhs_n$} \label{sec:lhsrhs}
From now on, unless specified, we consider trees to be labeled and rooted. 

A {\em tree} on $[n]:=\{1,2,\ldots,n\}$ is an acyclic connected graph on the vertex set~$[n]$ such that one vertex, called the {\em root}, is distinguished. We denote by~$\T_{n}$ the set of trees on $[n]$ and by $\T_{n,i}$ the set of trees on $[n]$ where vertex $i$ is the root. A~{\em forest} is a graph such that every connected component is a tree. Let $\For_{n}$ denote the set of forests on $[n]$. There is a canonical bijection $\gamma: \T_{n+1,n+1} \to \For_n$ such that $\gamma(T)$ is the forest obtained from $T$ by removing the vertex $n+1$ and letting each neighbor of $n+1$ be a root. A graph is called {\em bicolored} if each vertex is colored with the color $\one$ (black) or $\two$ (white). We denote by $\For^{\bi}_n$ the set of bicolored forests on $[n]$. From Cayley's formula~\cite{C} and the bijection $\gamma$, we have 
\begin{equation}
|\For_n|=|\T_{n+1,n+1}|=(n+1)^{n-1} \quad \mbox{and} \quad |\For^{\bi}_n|=2^{n}\cdot (n+1)^{n-1}.
\label{equ:cay}
\end{equation} 
Thus $\lhs_n$ can be interpreted as the cardinality of $\For_n^{\bi}$.

Let $F$ be a forest and let $i$ and $j$ be vertices of $F$. We say that $j$ is a {\em descendant} of $i$ if $i$ is contained in the path from $j$ to the root of the component containing $j$. In particular, if $i$ and $j$ are joined by an edge of $F$, then $j$ is called a {\em child} of $i$. Note that~$i$ is also a descendant of $i$ itself. Let $S(F,i)$ be the induced subtree of $F$ on descendants of $i$, rooted at $i$. We call this tree the descendant subtree of $F$ rooted at $i$. A vertex $i$ is called {\em proper} if $i$ is the smallest vertex in $S(F,i)$\,; otherwise $i$ is called {\em improper}. Let $\pv(F)$ denote the the number of proper vertices in $F$.

A {\em plane tree} or {\em ordered  tree} is a tree such that the children of each vertex are linearly ordered. We denote by $\Pla_{n}$ the set of plane trees on $[n]$ and by $\Pla_{n,i}$ the set of plane trees on $[n]$ where vertex $i$ is the root. Define a {\em plane forest} on $[n]$ to be a finite sequence of non-empty plane trees $(P_{1},\ldots,P_{m})$ such that~$[n]$ is the disjoint union of the sets $V\!(P_{r})$,~$1\leq r \leq m$\,. We denote by~$\PF_n$ the set of plane forests on $[n]$ and by $\PF^{\bi}_n$ the set of bicolored plane forests on $[n]$. There is also a canonical bijection ${\bar \gamma}: \Pla_{n+1,n+1}\to \PF_{n}$ such that ${\bar \gamma}(P)=\big(S(P,j_{1}),\ldots,S(P,j_{m})\big)$ where each vertex $j_{r}$ is the  $r$th child of $n+1$ in $P$. It is well-known that the number of unlabeled plane trees on $n+1$ vertices is given by the $n$th Catalan number $C_n=\frac{1}{n+1}\binom{2n}{n}$ (see~\cite[ex.~6.19]{St}). Thus we have 
\begin{equation}
|\PF_n|=|\Pla_{n+1,n+1}|=n!\cdot C_n = 2n\,(2n-1)\cdots(n+2)\,.
\label{equ:ord}
\end{equation}

A {\em binary tree} is a tree in which each vertex has at most two children and each child of a vertex is designated as its left or right child. We denote by $\Bin_n$ the set of binary trees on $[n]$ and by $\Bin^{\bi}_n$ the set of  bicolored binary trees on $[n]$. 

For $k \geq 2$\,, a {\em $k$-ary tree} is a tree where each vertex has at most $k$ children and each child of a vertex is designated as its first, second,\,\dots, or $k$th child. We denote by $\Ary^{k}_n$ the set of $k$-ary trees on $[n]$. Clearly, we have that $\Ary^{2}_{n} =\Bin_n$. Since the number of unlabeled $k$-ary trees on $n$ vertices is given by $\frac{1}{(k-1)n+1}\binom{kn}{n}$ (see \cite[p.~172]{St}), the cardinality of $\Ary^{k}_n$ is as follows:

\begin{equation}
|\Ary^{k}_n|=n!\cdot \frac{1}{(k-1)n+1}\binom{kn}{n} = kn\,(kn-1)\cdots(kn-n+2)\,.
\end{equation}

Now we introduce a combinatorial interpretation of the number $\rhs_n$. Let $\ubt$ be an unlabeled binary tree on $n$ vertices and $\omega:V\!(\ubt) \to [n]$ be a bijection. Then the pair $(\ubt, \omega)$ is identified with a (labeled) binary tree on $[n]$. Let $\Pi(\lbt)$ be the set of vertices $v$ in $\ubt$ such that $v$ has no descendant $v'$ satisfying $\omega(v)>\omega(v')$\,.

Let $\D_{n}$ be the set of bicolored binary trees on $[n]$ such that each proper vertex is colored with $\one$ or $\two$ and each improper vertex is colored with $\one$.

\begin{lemma} \label{lem:rhs}
The cardinality of $\D_n$ is equal to $\rhs_n$.
\end{lemma}
\begin{proof}
Let $\D'_{n}$ be the set defined as follows:
$$
\D'_{n}:=\{\,(\lbt,\alpha) \mid (\lbt) \in \Bin_n \mbox{\rm~and~} \alpha \subseteq \Pi(\lbt) \,\} \,.   
$$
There is a canonical bijection from $\D'_n$  to $\D_n$ as follows: Given $(\lbt,\alpha) \in \D'_n$\,, if a vertex $v$ of $\ubt$ is contained in $\alpha$ then color $v$ with $\two$; otherwise color $v$ with~$\one$. Thus it suffices to show that the cardinality of $\D'_n$ equals $\rhs_n$.
  
Given an unlabeled binary tree $\ubt$ and a subset $\alpha$ of $V(\ubt)$, let $l(\ubt, \alpha)$ be the number of labelings $\omega$ satisfying $\alpha \subseteq \Pi(\lbt)$\,. Then for each $v \in \alpha$\,, the label~$\omega(v)$ of $v$ should be the smallest label among the labels of the descendants of~$v$. So the number of possible labelings $\omega$ is $n!/ \prod_{v \in \alpha} h(v)$\,. Thus we  have
\begin{eqnarray*}
|\D'_{n}| 
&=& \sum_{\ubt} \sum_{\alpha \subseteq V\!(\ubt)} l(\,\ubt,\alpha)\\
&=& \sum_{\ubt} \sum_{\alpha \subseteq V\!(\ubt)} {n!}\,\prod_{v \in \alpha} \, \dfrac{1}{h(v)}\,, \\
\end{eqnarray*} 
which coincides with $\rhs_n$.   
\end{proof}

\section{A bijection}  \label{sec:bij}
In this section, we construct a bijection between $\For_n^{\bi}$ and $\D_n$, which gives a bijective proof of (\ref{eqn:main}). 

Given a vertex $v$ of a bicolored binary tree $B$, let $L(B,v)$ (resp.~$R(B,v)$\,) be the descendant subtree of $B$, which is rooted at the left (resp.~right) child of $v$. Note that $L(B,v)$ and $R(B,v)$ may be empty, but $L(B,v)$ or $R(B,v)$ is nonempty when $v$ is improper. For any kind of tree $T$, let $m(T)$ be the smallest vertex in~$T$. By convention, we put $m(\emptyset)=\infty$\,. For an improper vertex $v$ of $B$, if $m\big(L(B,v)\big)>m\big(R(B,v)\big)$\,, then we say that $v$ is {\em right improper}\,; otherwise {\em left improper}. 

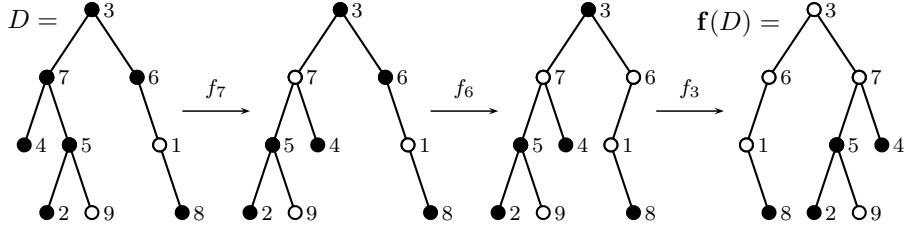
\begin{figure}
\centering
\begin{pspicture}(0,0)(12,3.6)\footnotesize
\psset{linewidth=1\pslinewidth}
\psset{dotsize=.2 0} 
\psset{xunit=.3cm, yunit=.9cm}
\rput(4,3.5){%
 \rput[tr](-1.5,0){{\normalsize $D=$}}
 \psline{*-*}(0,0)(-2,-1)
 \psline{*-*}(0,0)(2,-1)
 \psline{*-*}(-2,-1)(-3,-2)
 \psline{*-*}(-2,-1)(-1,-2)
 \psline{*-o}(2,-1)(3,-2)
 \psline{*-*}(-1,-2)(-2,-3)
 \psline{*-o}(-1,-2)(-0,-3)
 \psline{o-*}(3,-2)(4,-3)
 \rput(0.75,0){%
  \rput(0,0){$3$}
  \rput(-2,-1){$7$}
  \rput(2,-1){$6$}
  \rput(-3,-2){$4$}
  \rput(-1,-2){$5$}
  \rput(3,-2){$1$}
  \rput(-2,-3){$2$}
  \rput(0,-3){$9$}
  \rput(4,-3){$8$}
 }
}
\psline[linewidth=.7\pslinewidth]{->}(8,2)(11,2)
\rput[b](9.5,2.2){$f_7$}
\rput(15,3.5){%
 \psline{*-o}(0,0)(-2,-1)
 \psline{*-*}(0,0)(2,-1)
 \psline{o-*}(-2,-1)(-3,-2)
 \psline{o-*}(-2,-1)(-1,-2)
 \psline{*-o}(2,-1)(3,-2)
 \psline{*-*}(-3,-2)(-4,-3)
 \psline{*-o}(-3,-2)(-2,-3)
 \psline{o-*}(3,-2)(4,-3)
 \rput(0.75,0){%
  \rput(0,0){$3$}
  \rput(-2,-1){$7$}
  \rput(2,-1){$6$}
  \rput(-3,-2){$5$}
  \rput(-1,-2){$4$}
  \rput(3,-2){$1$}
  \rput(-4,-3){$2$}
  \rput(-2,-3){$9$}
  \rput(4,-3){$8$}
 }
}
\psline[linewidth=.7\pslinewidth]{->}(19,2)(22,2)
\rput[b](20.5,2.2){$f_6$}
\rput(26,3.5){%
 \psline{*-o}(0,0)(-2,-1)
 \psline{*-o}(0,0)(2,-1)
 \psline{o-*}(-2,-1)(-3,-2)
 \psline{o-*}(-2,-1)(-1,-2)
 \psline{o-o}(2,-1)(1,-2)
 \psline{*-*}(-3,-2)(-4,-3)
 \psline{*-o}(-3,-2)(-2,-3)
 \psline{o-*}(1,-2)(2,-3)
 \rput(0.75,0){%
  \rput(0,0){$3$}
  \rput(-2,-1){$7$}
  \rput(2,-1){$6$}
  \rput(-3,-2){$5$}
  \rput(-1,-2){$4$}
  \rput(1,-2){$1$}
  \rput(-4,-3){$2$}
  \rput(-2,-3){$9$}
  \rput(2,-3){$8$}
 }
}
\psline[linewidth=.7\pslinewidth]{->}(29,2)(32,2)
\rput[b](30.5,2.2){$f_3$}
\rput(36,3.5){%
 \rput[tr](-1.5,0){{\normalsize $\f(D)=$}}
 \psline{o-o}(0,0)(-2,-1)
 \psline{o-o}(0,0)(2,-1)
 \psline{o-o}(-2,-1)(-3,-2)
 \psline{o-*}(2,-1)(1,-2)
 \psline{o-*}(2,-1)(3,-2)
 \psline{o-*}(-3,-2)(-2,-3)
 \psline{*-*}(1,-2)(0,-3)
 \psline{*-o}(1,-2)(2,-3)
 \rput(0.75,0){%
  \rput(0,0){$3$}
  \rput(-2,-1){$6$}
  \rput(2,-1){$7$}
  \rput(-3,-2){$1$}
  \rput(1,-2){$5$}
  \rput(3,-2){$4$}
  \rput(-2,-3){$8$}
  \rput(0,-3){$2$}
  \rput(2,-3){$9$}
 }
}
\end{pspicture}
\caption{The flip $\f$ of $D$} \label{fig:flip}
\end{figure}

For a vertex $v$ of $B$, define the {\em flip} on $v$, which will be denoted by $f_v$\,, by swapping $L(B,v)$ and $R(B,v)$ and changing the color of $v$. Note that this flip operation satisfies $f_v \circ f_v =id$ and $f_{v}\circ f_{w} = f_{w}\circ f_{v}$\,. For a bicolored binary tree $D$ in $\D_n$, let $\f$ be the map defined by 
$$
\f(D):=(f_{v_{1}}\circ \cdots \circ f_{v_{k}})(D)\,, 
$$
where $\{v_{1},\ldots,v_{k}\}$ is the set of right improper vertices in $D$. (See Figure~\ref{fig:flip}.) 

Let $\E_{n}$ be the set of bicolored binary trees $E$ on~$[n]$ such that every improper vertex $v$ is left improper, i.e., $m\big(L(E,v)\big)<m\big(R(E,v)\big)$.
\begin{lemma}
The map $\f$ is a bijection from $\D_n$ to $\E_n$. 
\label{lem:DE}
\end{lemma}
\begin{proof}
For a bicolored binary tree $E$ in $\E_n$, let $\f'$ be the map defined by $\f'(E):=(f_{u_{1}}\circ \cdots \circ f_{u_{j}})(E)$\,, where $\{u_{1},\ldots,u_{j}\}$ is the set of white-colored improper vertices in $E$.

Let $v$ be a right (resp. white-colored) improper vertex of $B\in\Bin^{\bi}_n$. Then~$v$ is a left  (resp. black-colored) improper vertex of $f_v(B)$. Moreover, since~$f_v$ does not change the state of other vertices, $f_v(B)$ has one less right (resp. white-colored) improper vertex than $B$. Thus the set of right improper vertices in $D$ equals the set of white-colored improper vertices in $\f(D)$ and $\f(D)$ is contained in $\E_n$. Similarly, the set of white-colored improper vertices in $E$ equals the set of right improper vertices in $\f'(E)$ and $\f'(E)$ is contained in $\D_n$. Since flip operations are commutative, we have that $(\f'\circ \f)(D)=D$ and  $(\f\circ \f')(E)=E$ for all $D\in \D_n$ and $E\in \E_n$, which completes the proof. 
\end{proof}

Let $\G_n$ (resp. $\Q_n$) be the set of bicolored trees (resp. bicolored plane trees) on~$[n+1]$ such that $n+1$ is the root colored with $\one$. Note that the map~$\gamma$~(resp.~$\bar{\gamma}$) can be regarded as a bijection $\gamma:\G_n \to \For_n^{\bi}$~(resp.~$\bar{\gamma}: \Q_n \to \PF_n^{\bi}$). It is easy to show that $\G_n$ can be viewed as a subset of $\Q_n$ satisfying the following condition: For an interior vertex $v$ of~$Q\in\Q_n$\,, let $(w_{1},\ldots, w_{r})$ be the children of $v$, in order. Then  $m\big(S(Q,w_{1})\big)<\cdots<m\big(\,S(Q,w_{r})\big)$ holds.

\begin{figure}
\centering
\begin{pspicture}(0,0)(12,4)\footnotesize
\psset{linewidth=1\pslinewidth}
\psset{dotsize=.2 0} 
\psset{xunit=.6cm, yunit=1cm}
\rput(4,3.5){%
 \rput[r](-1.5,0){{\normalsize $E=$}}
 \psline{o-o}(0,0)(-2,-1)
 \psline{o-o}(0,0)(2,-1)
 \psline{o-o}(-2,-1)(-3,-2)
 \psline{o-*}(2,-1)(1,-2)
 \psline{o-*}(2,-1)(3,-2)
 \psline{o-*}(-3,-2)(-2,-3)
 \psline{*-*}(1,-2)(0,-3)
 \psline{*-o}(1,-2)(2,-3)
 \rput(0.6,0){%
  \rput(0,0){$3$}
  \rput(-2,-1){$6$}
  \rput(2,-1){$7$}
  \rput(-3,-2){$1$}
  \rput(1,-2){$5$}
  \rput(3,-2){$4$}
  \rput(-2,-3){$8$}
  \rput(0,-3){$2$}
  \rput(2,-3){$9$}
 }
}
\psline[linewidth=.7\pslinewidth]{->}(8.4,2)(11,2)
\rput[b](9.7,2.2){{\normalsize $\phi$}}
\rput(15.5,3.5){%
 \rput[r](-1.5,0){{\normalsize $\phi(E)=$}}
 \psline{*-o}(0,0)(-3,-1)
 \psline{*-o}(0,0)(0,-1)
 \psline{*-*}(0,0)(3,-1)
 \psline{o-o}(-3,-1)(-3,-2)
 \psline{o-*}(0,-1)(-1,-2)
 \psline{o-o}(0,-1)(1,-2)
 \psline{o-o}(-3,-2)(-4,-3)
 \psline{o-*}(-3,-2)(-2,-3)
 \psline{*-*}(-1,-2)(-1,-3)
  \rput(0.6,0){%
  \rput(0.1,0){$10$}
  \rput(-3,-1){$3$}
  \rput(0,-1){$7$}
  \rput(3,-1){$4$}
  \rput(-3,-2){$6$}
  \rput(-1,-2){$5$}
  \rput(1,-2){$9$}
  \rput(-4.1,-3){$1$}
  \rput(-3.1,-3){$8$}
  \rput(-1,-3){$2$}
 }
}
\psline[linewidth=.7\pslinewidth]{->}(29,2)(32,2)
\rput[b](30.5,2.2){$f_3$}
\end{pspicture}
\caption{The bijection $\phi$} \label{fig:phi}
\end{figure}
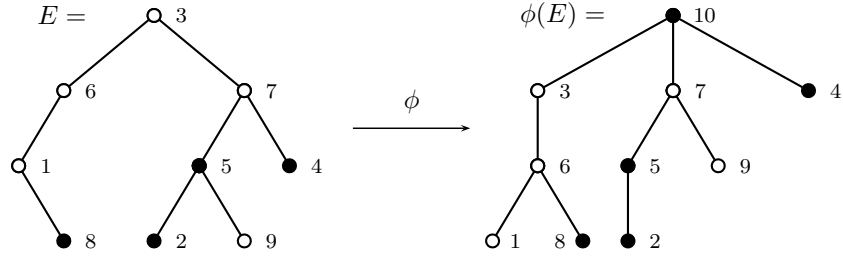

Recall that $\Bin^{\bi}_n$ denotes the set of bicolored binary trees on~$[n]$. Clearly we have~$\E_n \subseteq \Bin^{\bi}_n$ and $\G_n \subseteq \Q_n$\,. Let $\Phi$ be a bijection from $\Bin^{\bi}_n$ to~$\Q_n$, which maps~$B$ to $Q$ as follows:
\begin{enumerate}
\item The vertices of $B$ are the vertices of $Q$ with the root deleted.
\item The root of $B$ is the first child of the black root $n+1$ of $Q$.
\item $v$ is a left child of $u$ in $B$ iff $v$ is the first child of $u$ in $Q$.
\item $v$ is a right child of $u$ in $B$ iff $v$ is the sibling to the right of $u$ in $Q$. 
\item The color of $v$ in $B$ is the same as the color of $v$ in $Q$.
\end{enumerate}
Note that here $\Phi$ is essentially an extension of a well-known bijection, which is described in~\cite[p.~60]{SW}, from binary trees to plane trees. 
\begin{lemma} \label{lem:EC}
The restriction $\phi$ of $\Phi$ to $\E_n$ is a bijection from $\E_n$ to $\G_n$. 
\end{lemma} 
\begin{proof}
For any improper vertex $v$ of $E \in \E_n$, we have $m(L(E,v))<m(R(E,v))$.  This guarantees that $m(S(G,v))<m(S(G,w))$ in $G=\phi(E)$, where $w$ (if it exists) is the sibling to the right of $v$ in $G$. Thus $\phi(E) \in \G_n$, i.e., $\phi(\E_{n})\subseteq \G_n$. Similarly we can show that $\phi^{-1} (\G_{n}) \subseteq \E_n$. So we have $\phi(\E_{n})= \G_n$, which implies that  $\phi$ is bijective. (See Figure~\ref{fig:phi}.)
\end{proof}
From Lemma~\ref{lem:EC}, we easily get that $\gamma\circ\phi$ is a bijection from $\E_n$ to $\For^{\bi}_n$. Combining this result with Lemma~\ref{lem:DE} yields the following consequence.
\begin{theorem} 
The map $\gamma\circ \phi \circ \f$ is a bijection from $\D_n$ to $\For^{\bi}_n$. \label{thm:bij}
\end{theorem} 
From equation (\ref{equ:cay}) the cardinality of $\For^{\bi}_n$ equals $\lhs_n$ and from Lemma~\ref{lem:rhs} the cardinality of $\D_n$ equals $\rhs_n$. Thus Theorem \ref{thm:bij} is a combinatorial explanation of identity (\ref{eqn:main}).

\section{Generalized formulas} \label{sec:gen}
In Theorem~\ref{thm:bij}, we showed that the set $\D_n$ of binary trees on~$[n]$ such that each proper vertex is colored with the color $\one$ or $\two$ and each improper vertex is colored with the color $\one$ has cardinality $|\D_n|=2^n\,(n+1)^{n-1}$\,. In this section, we give a generalization of this result. 

For $n \geq 1$, let $a_{n,m}$ denote the number of $k$-ary trees on~$[n]$ with~$m$ proper vertices.  By convention, we put $a_{0,m}=\delta_{0,m}$\,.  Let
$$
a_n(t)=\sum_{m\geq 0} a_{n,m}\, t^{m}=\sum_{T \in \Ary^k_n} t^{\,\pv(T)}\,.
$$
It is clear that for a positive integer $t$ the number $a_n(t)$ is the number of $k$-ary trees on~$[n]$ such that each proper vertex is colored with the color $\1$, $\2$,\,\dots, or $\bt$ and each improper vertex has  one color $\1$. 
Let $A(x)$ be denote the exponential generating function for $a_n(t)$, i.e.,
$$
A(x)=\sum_{n\geq 0} a_n(t)\, \frac{x^n}{n!}\,.
$$
\begin{lemma}
The function $A=A(x)$ satisfies the following equation:
\begin{equation}
A=\big(\,1+(kt-t-k)\,x\,A^{k-1}\,\big)^{t/(kt-t-k)}.
\label{funeq:ary}
\end{equation}
\label{lem:funeq}
\end{lemma}
\begin{proof}
Let $T$ be an $k$-ary tree on $[n]\cup \{0\}$. Delete all edges going from the root~$r$ of $T$. Then $T$ is decomposed into $T'=(r;T_{1},\ldots,T_{k})$ where each $T_{i}$ is a $k$-ary tree and $[n]\cup \{0\}$ is the disjoint union of $V(T_{1}),\ldots,V(T_{k})$ and $\{r\}$. Consider two cases:
(i) For some $1\leq i \leq k$, $T_{i}$ has the vertex $0$\,; (ii) $r=0$.  Then we have
\begin{eqnarray*}
a_{n+1}(t)
&\!\!=\!\!&\sum_{i=1}^{k}~\sum_{n_{1}+\cdots+n_{k}=n-1}\binom{n}{1,n_{1},\ldots,n_{k}} \,a_{n_1}(t)\,\cdots\, a_{n_{i}+1}(t)\,\cdots\, a_{n_k}(t) \\
&& +\quad t \sum_{n_{1}+\cdots+n_{k}=n} \binom{n}{n_{1},\ldots,n_{k}}\,a_{n_1}(t)\,\cdots \,a_{n_k}(t) \,.
\end{eqnarray*}
Multiplying both sides by $x^{n}/n!$ and summing over $n$ yields
\begin{equation} \label{diff:ary}
A' = k\,x\,A^{k-1}A' + t\,A^{k} \,, 
\end{equation}
with $A(0)=1$, where the prime denotes the derivative with respect to $x$.
Adding~$(kt-t-k)\,x\,A^{k-1}A'$ to both sides of (\ref{diff:ary}) yields
$$
\big(\,1+(kt-t-k)\,x\,A^{k-1}\,\big)\,A'=t\,\big(\,(k-1)\,x\,A^{k-2}A' + A^{k-1}\,\big)\,A\,.
$$ 
Let $\alpha(x)=1+(kt-t-k)\,x\,A^{k-1}$ and $\beta(x)=(k-1)\,x\,A^{k-2}A' + A^{k-1}$\,. Since~$\alpha'(x)=(kt-t-k)\,\beta(x)$, we have 
$$
(\,\log\,A\,)'=\frac{t}{kt-t-k}\,\big(\log\,\alpha(x)\,\big)',
$$
which implies the functional equation (\ref{funeq:ary}).
\end{proof}

Now we can deduce a formula for the polynomial~$a_n(k)$ from equation~(\ref{funeq:ary}).
\begin{theorem}[$k$-ary trees]
For $n\geq 1$\,, the polynomial $a_n(t)$ in $t$ is given by
\begin{equation}
a_n(t)=t\, \prod_{i=1}^{n-1}\big(\,(ki-i+1)\,t + k\,(n-i)\,\big)\,.
\label{equ:gen}
\end{equation}
\end{theorem}
\begin{proof}
Let $y^k=x$ and ${\sA}(y)=y\,A(y^{k-1})$. Then by Lemma~\ref{lem:funeq} we have 
\begin{equation} \label{li:ary}
{\sA}(y)=y\,\big(\,1+(kt-t-k)\,{\sA}(y)^{k-1}\,\big)^{t/(kt-t-k)}.
\end{equation}
Note that
\begin{equation*}
a_n(t)=\left[\dfrac{x^n}{n!}\right]A(x)=\left[\dfrac{y^{(k-1)n}}{n!}\right]A(y^{k-1})=n!\left[y^{kn-n+1}\right]{\sA}(y)\,.
\label{equ:coeff}
\end{equation*}
Applying the Lagrange Inversion Formula (see \cite[p.~38]{St}) to (\ref{li:ary}) yields that
\begin{eqnarray*}
\left[y^{kn-n+1}\right]{\sA}(y)
&=&\frac{1}{kn-n+1} \left[y^{(k-1)n}\right]
\big(\,1+(kt-t-k)\,y^{k-1}\,\big)^{\mbox{\small{$\frac{t\,(kn-n+1)} {kt-t-k}$}}}\\
&=&\frac{1}{kn-n+1}\,(kt-t-k)^{n}\,\binom{\frac{t\,(kn-n+1)}{kt-t-k}} {n}\\
&=&\frac{t}{n!}\,\prod_{i=1}^{n-1}\,\big(\,t\,(kn-k+1) - (kt-t-k)\,i\,\big)\,.
\end{eqnarray*}
Thus we obtain the desired result.
\end{proof}
Clearly, substituting $t=1$ in (\ref{equ:gen}) yields the number of $k$-ary trees on $[n]$, i.e.,
$$
a_n(1)=kn\,(kn-1)\cdots(kn-n+2)=|\Ary_n^k|\,.
$$ 
For some values of $k$, we can get interesting results. In particular when~$k=2$ we have   
$$ 
a_n(t)= t\prod_{i=1}^{n-1} \big(\,(i+1)t+2(n-i)\,\big)
~\stackrel{t=2} \longrightarrow~2^{n}(n+1)^{n-1}\,,
$$
which is a generalization of $|\D_n|=2^n\,(n+1)^{n-1}$, i.e., identity (\ref{eqn:main}).

For $n \geq 0$, let $f_{n,m}$ denote the number of forests on~$[n]$ with~$m$ proper vertices  and $p_{n,m}$ denote the number of plane forests on $[n]$  with~$m$ proper vertices . Let
$$
f_n(t)=\sum_{m\geq 0} f_{n,m}\, t^{m} 
\quad \mbox{~and~} \quad
p_n(t)=\sum_{m\geq 0} p_{n,m}\, t^{m}\,. 
$$
Let $F(x)$ and $P(x)$ be the exponential generating function for $f_n(t)$ and $p_n(t)$, respectively, i.e., 
$$
F(x)=\sum_{n\geq 0} f_n(t)\, \frac{x^n}{n!}
\quad \mbox{~and~}\quad
P(x)=\sum_{n\geq 0} p_n(t)\, \frac{x^n}{n!}\,.
$$
With the same methods used for~$k$-ary trees, we can get the following results.
\begin{theorem}[forests]
Suppose $f_n(t)$ and $F(x)$ are defined as above. Then we have
\begin{enumerate}
\item $F=F(x)$ satisfies the following differential equation:
$$F'=x\,F\,F' + t\,F^2\,, \quad \mbox{with~~$F(0)=1$}\,.$$
\item $F$ satisfies the following functional equation:
$$F=\big(\,1+(t-1)\,x\,F\,\big)^{t/(t-1)}\,.$$
\item For $n\geq 1$, the polynomial~$f_n(t)$ in~$t$ is given by
\begin{equation}
f_n(t)=t\, \prod_{i=1}^{n-1}\big(\,(i+1)\,t + (n-i)\,\big)\,.
\label{equ:gencay}
\end{equation}
\end{enumerate}
\end{theorem}
\begin{theorem}[plane forests]
Suppose $p_n(t)$ and $P(x)$ are defined as above. Then we have
\begin{enumerate}
\item $P=P(x)$ satisfies the following differential equation:
$$P'=x\,P^2 P' + t\,P^3\,, \quad \mbox{with~~$P(0)=1$}\,.$$
\item $P$ satisfies the following functional equation:
$$P=\big(\,1+(2t-1)\,x\,P^2\,\big)^{t/(2t-1)}\,.$$
\item For $n\geq 1$, the polynomial~$p_n(t)$ in~$t$ is given by
\begin{equation}
p_n(t)=t\, \prod_{i=1}^{n-1}\big(\,(2i+1)\,t + (n-i)\,\big)\,.
\label{equ:genord}
\end{equation}
\end{enumerate}
\end{theorem}
Note that the polynomials (\ref{equ:gencay}) and (\ref{equ:genord}) are generalizations of (\ref{equ:cay}) and (\ref{equ:ord}), respectively. Moreover, from these formulas, we can easily get
\begin{eqnarray*}
\sum_{T \in \T_{n+1}} t^{\,\pv(T)} &=& t\,\prod_{i=0}^{n-1}\big(\,(i+1)\,t + (n-i)\,\big), \\
\sum_{P \in \Pla_{n+1}} t^{\,\pv(P)} &=& t\,\prod_{i=0}^{n-1}\big(\,(2i+1)\,t + (n-i)\,\big),
\end{eqnarray*}
which are generalizations of $|\T_{n+1}|=(n+1)^n$ and $|\Pla_{n+1}|=(n+1)!\,C_n$\,.
\section*{Acknowledgment}
The author thanks Ira M. Gessel for his helpful advice and suggestions.


\end{document}